\documentclass{amsart}

\usepackage{amsmath,amssymb,amsthm}
\usepackage{graphicx,tikz}

\newtheorem*{thm}{Theorem}

\newtheorem*{corollary}{Corollary}

\theoremstyle{definition}

\theoremstyle{remark}

\DeclareMathOperator{\lcm}{lcm}

\begin{document}

\title[]{Intrinsic Sparsity of Kantorovich Solutions }
\thanks{The authors gratefully acknowledge support from the Kantorovich Initiative. S.S. was partially supported by the NSF (DMS-2123224) and the Alfred P. Sloan Foundation.}

\author[]{Bamdad Hosseini}
\address{Department of Applied Mathematics, University of Washington, Seattle, WA 98195}
\email{bamdadh@uw.edu }

\author[]{Stefan Steinerberger}
\address{Department of Mathematics, University of Washington, Seattle, WA 98195}
\email{steinerb@uw.edu}

\begin{abstract} Let $X,Y$ be two finite sets of points having $\#X = m$ and $\#Y = n$ points with $\mu = (1/m) \sum_{i=1}^{m} \delta_{x_i}$ and
$\nu = (1/n) \sum_{j=1}^{n} \delta_{y_j}$ being the associated uniform probability measures. 
A result of Birkhoff implies that if $m = n$, then the Kantorovich problem has a solution which also solves the Monge problem: optimal transport can be realized with a bijection $\pi: X \rightarrow Y$.  This is impossible when $m \neq n$.  We observe that when $m \neq n$, there exists a solution of the Kantorovich problem such that the mass of each point in $X$ is moved to at most $n/\gcd(m,n)$ different points in $Y$ and that, conversely, each point in $Y$ receives mass from at most $m/\gcd(m,n)$ points in $X$. 
\end{abstract}

\maketitle

\vspace{20pt}

 Let $\mu$ and $\nu$ be two (probability) measures. A classical question, due to Monge, is to understand the optimal way of mapping $\mu$ to $\nu$. If we denote the cost of transporting mass from $x$ to $y$ by $c(x,y)$, then the Monge problem asks for
$$ \inf_{T} \left\{ \int_X c(x, T(x)) d\mu(x): T_{*}(\mu) = \nu \right\} \qquad  \mbox{(Monge),}$$
where $T_{*}(\mu)$ denotes the push forward of $\mu$ by $T$.  This problem may not be solvable because such transport maps $T_{}$ may simply not exist. Kantorovich proposed to relax the problem and instead try to minimize
$$ \inf_{\gamma}  \int_{X \times Y} c(x, y) d\gamma(x,y) \qquad \qquad  \mbox{(Kantorovich),}$$
where $\gamma$ is a probability measure on $X \times Y$ having marginals $\mu$ and $\nu$.
There is a nice classical result linking these two problems in the discrete setting: if $\mu$ and $\nu$ are two uniform probability measures over two sets $X$ and $Y$ with $n$ elements, then it is known that these two problems coincide.
\begin{thm}[see e.g. \cite{brezis, merigot, peyre}] If $\mu = (1/n) \sum_{i=1}^{n} \delta_{x_i}$ and $\nu = (1/n) \sum_{i=1}^{n} \delta_{y_i}$, then there is a solution of the Kantorovich problem which also solves the Monge problem.
\end{thm}
The statement is independent of the transport costs $c(x_i, y_j)$.
The argument is as follows: the Kantorovich problem can, in the discrete setting, be formulated as a linear program over bistochastic matrices. A theorem of Birkhoff \cite{birk} (also attributed to K\"onig \cite{konig} and von Neumann \cite{von}) states that the bistochastic matrices are the convex hull of the permutation matrices. The minimum of a linear
program in a non-empty polyhedron is attained in an extremal point.\\

No such statement can be true when $m \neq n$: the two sets have different cardinalities and no bijection is possible. The goal of this short note is to point out that there nonetheless exists a particularly simple solution of the Kantorovich problem.

\begin{thm} Let $\mu = (1/m) \sum_{i=1}^{m} \delta_{x_i}$ and $\nu = (1/n) \sum_{i=1}^{n} \delta_{y_i}$. There is a solution of the Kantorovich problem such that mass from each point 
 in $X$ is moved to at most $n/\gcd(m,n)$ different points in $Y$ and that each point in $Y$ receives mass from at most $m/\gcd(m,n)$ points in $X$.
\end{thm}

Somewhat to our surprise, we were unable to find this simple but intriguing statement (illustrated in Fig. 1) in the literature. Besides its intrinsic appeal, it does seem like it could be
potentially useful insofar as it guarantees the existence of `sparse' solutions of the Kantorovich problem (with sparsity depending on $m,n$).

\vspace{-0pt}

\begin{center}
\begin{figure}[h!]
\begin{tikzpicture}
\node at (0,0) {\includegraphics[width=0.35\textwidth]{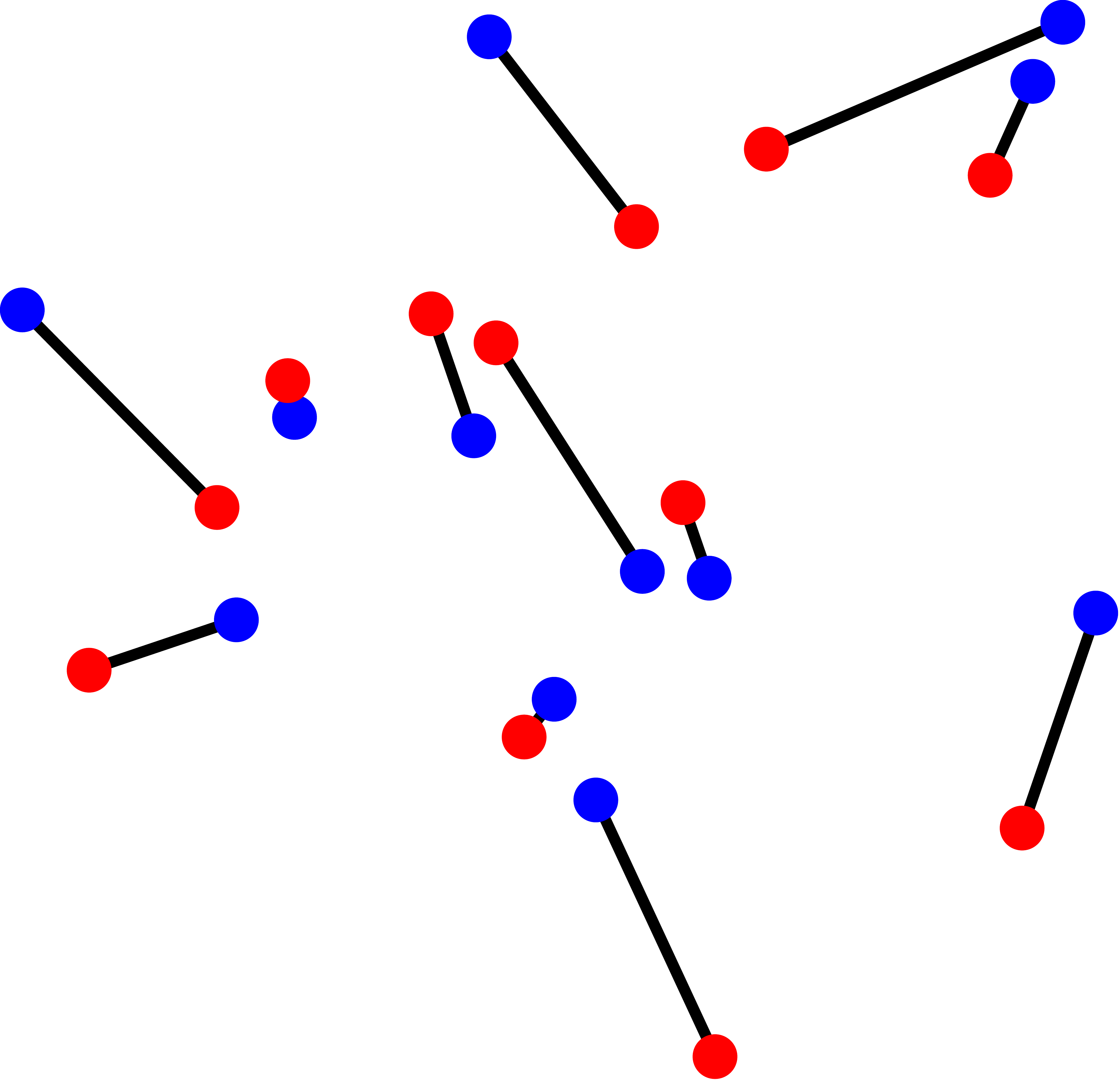}};
\node at (6,0) {\includegraphics[width=0.35\textwidth]{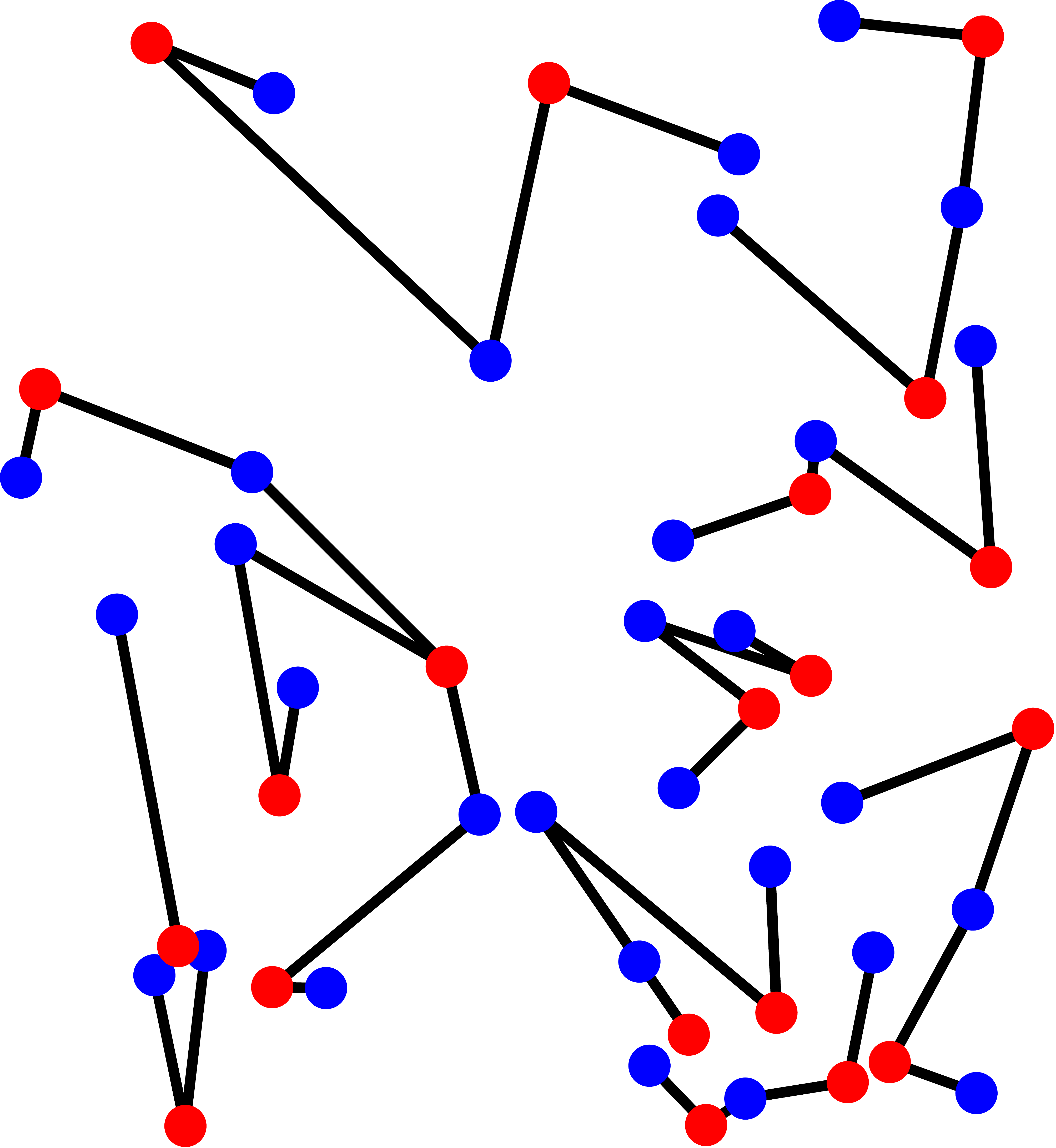}};
\end{tikzpicture}
\caption{Left: $m=n$, the transport is a bijection. Right: $m=20$ red points are sent to $n=30$ blue points. Each red point is transported to at most $30/\gcd(20,30) = 3$ blue points, each blue points receives mass from at most $20/\gcd(20,30) = 2$ red points.}
\end{figure}
\end{center}

\vspace{-20pt}

\begin{proof} 
Suppose that
$$ \mu = \frac{1}{m} \sum_{i=1}^{m} \delta_{x_i} \qquad \mbox{and} \qquad \nu = \frac{1}{n} \sum_{i=1}^{n} \delta_{y_i}$$
are two given measures. We replace each point $x_i$ by $n/\gcd(m,n)$ identical points $x^{}_{i,j}$ for $1\leq j \leq n/\gcd(m,n)$
and, likewise, we replace each point $y_i$ by $m/\gcd(m,n)$ identical points $y^{}_{i,j}$ where $1 \leq j \leq m/\gcd(m,n)$.
This allows us to write
$$ \mu = \frac{\gcd(m,n)}{mn} \sum_{i=1}^{m} \sum_{j=1}^{n/\gcd(m,n)} \delta_{x^{}_{i,j}} \qquad
\mbox{and} \qquad \nu = \frac{\gcd(m,n)}{mn} \sum_{i=1}^{n} \sum_{j=1}^{m/\gcd(m,n)} \delta_{y^{}_{i,j}}.$$
The problem can now be interpreted as finding a transport map from $mn/\gcd(m,n)$ points of the same weight to
another set of $mn/\gcd(m,n)$ points of the same weight. Applying the classical result shows that there exists
bijective map between the points that realizes the optimal Kantorovich cost. This corresponds into each point
in $X$ being split into at most $n/\gcd(m,n)$ equal parts and each point in $Y$ being split into at most $m/\gcd(m,n)$ parts. 
\end{proof}

We conclude by observing that the same argument also applies to linear combinations of weighted Dirac measures as long as the weights 
are rational. Suppose
$$ \mu = \sum_{i=1}^{m} \frac{a_i}{b_i} \delta_{x_i} \qquad \mbox{where} \qquad  \sum_{i=1}^{m} \frac{a_i}{b_i} = 1$$
and where $a_i,b_i \in \mathbb{N}$ are positive rational weights. We note that it is possible to equivalently represent $\mu$
as a linear combination of a number of equally weighted Dirac measures. This number will depend on the least common multiple  $\lcm(b_1, \dots, b_m)$
of the denominators. This can be seen by writing
$$ \mu = \frac{1}{\lcm(b_1, \dots, b_m)} \sum_{i=1}^{m} \frac{a_i \lcm(b_1, \dots, b_m)}{b_i}\delta_{x_i} $$
and noting that $a_i \lcm(b_1, \dots, b_n)/b_i \in \mathbb{N}$. This implies the following corollary.

\begin{corollary} Let 
$$\mu =  \sum_{i=1}^{m} \frac{a_i}{b_i} \delta_{x_i} \qquad \mbox{and} \qquad \nu =  \sum_{i=1}^{n} \frac{c_i}{d_i} \delta_{y_i}$$
be two probability measures with positive rational weights and let
$$ B = \lcm(b_1, \dots, b_m) \qquad \mbox{and} \qquad D = \lcm(d_1, \dots, d_n).$$
There exists a solution of the Kantorovich problem such
that mass from each point in $X$ is moved to at most $D/\gcd(B,D)$ different points in $Y$ and each point in $Y$ receives mass
from at most $B/\gcd(B,D)$ different points in $X$.
\end{corollary}
We note that this reduces to the previous Theorem when $a_i = c_i = 1$, $b_i = m$ and $d_i = n$. We also observe that it is quite
possible that $D/\gcd(B,D) \gg n$ and that $B/\gcd(B,D) \gg m$ (this will usually happen when the $b_i, d_i$ have many different prime
factors). In such a case, the statement would not say anything of interest.

\end{document}